%&tex
\normalbaselineskip=1.6\normalbaselineskip\normalbaselines
\magnification=1200

\def\pmb#1{\setbox0=\hbox{#1}%
 \kern-.025em\copy0\kern-\wd0
 \kern.05em\copy0\kern-\wd0
 \kern-.025em\raise.0433em\box0 }

\def\Z{{\bf Z}}

\def \bs{\bigskip}

\centerline{\bf Torsion points on $y^2=x^6+1$ }
\medskip
\centerline{\bf Jos\'e Felipe Voloch}
\bs

Let $C$ be the curve $y^2=x^6+1$ of genus 2 over a field of
characteristic zero. Consider $C$ embedded in its Jacobian $J$
by sending one of the points at infinity on $C$ to the origin of
$J$. In this brief note we show that the points of $C$ whose 
image on $J$ are torsion are precisely the two points at infinity, the
two points with $x=0$
and the six points with $y=0$. The finiteness of this set follows
from the Manin-Mumford conjecture proved by Raynaud [R]. Bounds for
this set follows from the work of Coleman [C] and Buium [B]. We will
follow Buium's approach enhanced by some calculations from [VW].
For the determination of the full torsion on other curves of genus 2
by a different method, see [BG].

The curve $C$ maps to $E^2$ 
where $E$ is $y^2=x^3+1$ by $\phi:(x,y) \mapsto ((x^2,y),(x^{-2},yx^{-3}))$.
This map factors through the embedding of $C$ in its Jacobian $J$, and
since $\phi$ maps the points at infinity on $C$ to torsion points on $E^2$,
it is enough to determine the points $P \in C$ with $\phi(P)$ torsion,
which is what we will do.
We will work over $\Z_7$ and compute the unramified torsion in $E^2$ that
lands in $C$. By a result of Coleman ([C]), this is enough for our
purposes (see also [B]).

The elliptic curve $E$ is a canonical lift of its reduction $E_0$ modulo 7
and, modulo 49, the unramified torsion in $E$ is
the image of the elliptic Teichm\"uller map $\tau:E_0({\bar{\bf F}}_7) 
\to E(W_2({\bar{\bf F}}_7))$ (see [VW]). Moreover, if $P \in E_0, P \ne 0$ 
then $\tau(P)$ has $x$-coordinate $(x,4x^{10}+x^7+2x^4+5x)$, where $x$ is the 
$x$-coordinate of $P$. Let $f(x)=4x^{10}+x^7+2x^4+5x$.
This statement follows from the proof of Proposition 4.2 of
[VW], there it is shown that the second Witt coordinate $x_1$of the 
$x$-coordinate of $\tau(P)$ satisfies $x_1'= f'(x)$ and has degree
at most $10$ in $x$. By looking at the
two-torsion, we get $x_1(-1)=x_1(-2)=x_1(-4)=0$, hence $x_1=f$.

Let $U$ be the affine open subset of $C$ where $x \ne 0, \infty$.
Note that if $(P,Q)$ is on the image of $U$ in $E^2$ 
then the product of the $x$-coordinates of $P$ and $Q$ is 1. If both 
$P$ and $Q$ are the elliptic Teichm\"uller lifts of their reduction modulo 7
and $(P,Q)$ is on $U$, we get
$(x^2,f(x^2))(x^{-2},f(x^{-2}))=1$ (product of Witt vectors of lenght two)
which gives $x^{14}f(x^{-2})+x^{-14}f(x^2)=0$. However,
$x^{14}f(x^{-2})+x^{-14}f(x^2)= (x^6+1)^4/x^{12}$, hence the torsion points
on $U$ are precisely the six points with $y=0$, hence the result.

Of course, the above calculations use extensively the special features
of the curve in question. Another example where this technique can be
employed is to compute the points on $X: x^4+y^4=1$ which map to 
torsion points on $F^2$, where $F$ is the
elliptic curve $y^2=1-x^4$ and the map is $(x,y) \mapsto ((x,y^2),(x^2,y))$.
Note that $F^2$ is just a quotient of the Jacobian of $X$, which is 
3-dimensional. We work 5-adically and again use the results of [C] (noticing
that $F$ is ordinary at 5) to reduce to unramified points. 
The elliptic curve $F$ is a canonical lift of its reduction $F_0$ modulo 5
and the elliptic Teichm\"uller map has $x$-coordinate $(x,2x^9-2x)$.
If a point $(P,Q)$ on $F^2$ is on the image of $X$, the $x$-coordinate of $Q$
is the $x$-coordinate of $P$ squared. If both
$P$ and $Q$ are the elliptic Teichm\"uller lifts of their reduction modulo 5
we get $(x,2x^9-2x)^2 = (u,2u^9-2u)$ as Witt vectors of lenght two. Hence
$u=x^2$ and $2x^5(2x^9-2x) = 2u^9-2u$, which implies $x^4=\pm 1$, hence
these points and the points at infinity are the points on $X$ that map to
torsion points on $F^2$. For torsion points on Fermat curves embedded in
their Jacobians, see [CTT].

{\bf Acknowledgements:}
The author would like to thank
the NSA (grant MDA904-97-1-0037) for financial support.

\bigskip
\centerline{\bf References.}
\bigskip
\noindent
[B] A. Buium {\it Geometry of $p$-jets}, Duke Math. Jour. {\bf 82}
(1996), 349--367.
\medskip
\noindent
[BG] J. Boxall and D. Grant, {\it Examples of torsion points on genus 
two curves}, preprint, 1997.
\medskip
\noindent
[C] R. F. Coleman {\it Ramified torsion points on curves } Duke
Math. J. {\bf 54} (1987) 615--640.
\medskip
\noindent
[CTT] R. F. Coleman, A. Tamagawa and P. Tzermias {\it The cuspidal
torsion packet on the Fermat curve}, preprint, 1997.
\medskip
\noindent
[R] M. Raynaud {\it Courbes sur une vari\'et\'e ab\'elienne et points
de torsion} Invent. Math. {\bf 71}(1983)207--233.
\medskip
\noindent
[VW] J. F. Voloch and J. L. Walker, {\it Euclidean weights of codes 
from elliptic curves over rings}, preprint, 1997.
\medskip
\noindent

Dept. of Mathematics, Univ. of Texas, Austin, TX 78712, USA
\smallskip
\noindent
e-mail: voloch@math.utexas.edu

\end